\documentstyle[12pt,oneside]{amsart}

\title{Jet schemes of homogeneous hypersurfaces}
\author{Shihoko Ishii}
\author{Akiyoshi Sannai}
\author{Kei-ichi Watanabe} 
\address{Department of Mathematical Science, University of  Tokyo, 
Komaba, Meguro, Tokyo, 153-8914, Japan,
\newline
e-mail :  shihoko@@ms.u-tokyo.ac.jp}
\address{Deparment of Mathematical Science,  University of Tokyo, Komaba, Meguro,
Tokyo, 153-8914, Japan,
\newline
e-mail :  sannai@@ms.u-tokyo.ac.jp}
\address{Department of Mathematics, College of Human and Science, Nihon University, Setagaya, Tokyo, 156-0045, Japan,
\newline
e-mail : watanabe@@math.chs.nihon-u.ac.jp}

\newcommand{\bQ}{{\Bbb Q}}

\newcommand{\bN}{{\Bbb N}}
\newcommand{\bA}{{\Bbb A}}
\newcommand{\fa}{{\frak a}}

\newcommand{\bx}{{\bf x}}
\newcommand{\bz}{{\bf z}}
\newcommand{\bo}{{\bf 0}}
\newcommand{\fm}{{\frak m}}
\newcommand{\codim}{\operatorname{codim}}

\newcommand{\ord}{\operatorname{ord}}

\let \cedilla =\c
\renewcommand{\c}[0]{{\mathbb C}}  

\renewcommand{\o}[0]{{\mathcal O}} 

\newcommand{\q}[0]{{\mathbb Q}}

\newcommand{\spec}[0]{\operatorname{Spec}}

\def\to {\longrightarrow}

\newtheorem{thm}{Theorem}[section]

\newtheorem{lem}[thm]{Lemma}
\newtheorem{cor}[thm]{Corollary}

\newtheorem{problem}[thm]{Problem}

\theoremstyle{definition}
\newtheorem{defn}[thm]{Definition}

\newtheorem{say}[thm]{}

\newtheorem{rem}[thm]{Remark}

\theoremstyle{remark}

\begin{document}
\maketitle
\footnote{S. Ishii: partially supported by JSPS Grant-In-Aid (B)22340004, (S)19104001\\
A. Sannnai: partially supported by JSPS research fellow 08J08285,\\
K-i. Watanabe: partially supported by JSPS Grant-in-Aid (C) 20540050
 }

\begin{abstract}
  This paper studies the singularities of jet schemes of homogeneous hypersurfaces of general type.
  We obtain the condition of the degree and the dimension for the singularities of the jet schemes to be of dense $F$-regular type.
  This provides us with examples of singular varieties whose $m$-jet schemes have rational singularities for every $m$. 
\end{abstract}

\section{Introduction}

\noindent
  The concept  jet schemes   over an algebraic variety 
  was introduced by Nash 
 in his preprint in 1968 
 which is later  published as \cite{nash}.
   These spaces  
  represent
  the nature of the singularities of the base space. 
 In fact, papers \cite{ein}, \cite{e-Mus}, \cite{must01}, \cite{must02} 
  by  Musta\cedilla{t}\v{a}, Ein and Yasuda show that  geometric 
  properties of the jet schemes determine  properties 
  of the singularities 
  of the base space. 
  To summarize, their results among others are as follows:

{\it Let \( X \) be a variety of locally  a complete intersection  over an 
 algebraically closed field of characteristic zero. 
Then \( X_{m} \) is of pure dimension (resp. irreducible, normal) for all \( m\geq 1  \)
 if and only if \( X  \) has log-canonical (resp. canonical, terminal) singularities.}

Accoding to this form, it is natural to formulate the question: 
\begin{problem}
Does the following hold?
$X$ is non-singular if and only if $X_m$ has at worst certain ``mild" singularities for every $m\geq 1$.

Does the bound of a certain invariant of the singularities on $X_m$ 
characterize the somoothness of $X$?
\end{problem}
The easiest candidate for ``certain mild singularities" is a rational singularity.
  In this paper, we show that a rationality is not appropriate for the required statement in the problem.
  This is proved by providing with counter examples. 
  We  study the singularities of the jet schemes of homogeneous hypersurface of general type and obtain the condition of the degree and the dimension for the singularities of the jet schemes to be dense $F$-regular type.
  \begin{thm}
  Let $X$ be a  hypersurface in $\bA_k^N$ over a field $k$ of characteristic $0$ defined by a  polynomial of general type of 
  degree $d$,
  i.e., 
  $$f(x_1, \ldots, x_N)=\sum_{i_1\leq i_2,\leq \ldots\leq i_d}\xi_{i_1,i_2,\ldots, i_d}x_{i_1}x_{i_2}\cdots x_{i_d},$$
  where $\{\xi_{i_1,i_2,\ldots, i_d}\}$ are algebraically independent over $\bQ$.
  If $d^2\leq N$, then the jet scheme $X_m$ has at worst rational singularities for every 
  $m\in \bN$.
  \end{thm}

 A rational singularity is defined by using a resolution of the singularities. Since it is almost impossible to construct a resolution of the singularities
of the jet scheme even for the simplest singularities on the base variety because of too many variables on the jet scheme,
 we use the positive characteristic method. 
  The theorem shows  examples of singular $X$ whose jet schemes $X_m$ for all $m$ have  at worst rational singularities.
  We also show that $X$ is non-singular if and only if the $F$-pure threshold 
  does not change between $X_m$'s  for different $m$:
  
  \begin{thm}
  Let $X$ be a variety of locally a complete intersection at $0$ over a field of characteristic $p>0$.
  For $m<m'$, assume also $X_m, X_{m+1},\ldots, X_{m'}$ are  complete intersections at the trivial jets \\ $0_m,\ldots,  0_{m'}$. 
  Then, the following are equivalent:
\begin{enumerate}
\item[(i)]  
  $(X,0)$ is non-singular;
\item[(ii)]
   $fpt(X_m,0_m,0_m)=
  fpt(X_{m'},{\psi_{m'm}}^{-1}(0_m), 0_{m'})$, where $\psi_{m'm}:X_{m'}\to X_m$ is the truncation morphism.
\end{enumerate}
\end{thm}

Throughout this paper the base field $k$ is an algebraically closed field.  
  
\section{Preliminaries on jet schemes and positive characteristic methods}

\begin{say}
\label{notation}
For a scheme $X$ of finite type over an algebraically closed field $k$, 
we can associate the space of $m$-jet (or the $m$-jet scheme) $X_m$ for every  $m\in \bN$.  
The exact definition of the $m$-jet scheme and the basic properties can be seen in \cite{cr}.
We use the notation and the terminologies in \cite{cr}.
The canonical projection $X_m\to X$ is denoted by $\pi_m$.

If $X$ is a closed subscheme of  $\bA_k^N$ defined by an equation $f=0$, then the $m$-jet scheme $X_m$ is defined in $$\bA_k^{(m+1)N}=\spec k[x_i^{(j)}\mid 1\leq i\leq N, j=0,1,\ldots m]$$ by the equations 
$\{ F^{(j)}=0\}_{j=0,\ldots, m}$.
Here, the $F^{(j)}\in   k[x_i^{(j)}\mid 1~\leq ~i~\leq~N, j=0,1,\ldots m]$ is defined as follows:

$$f(\sum_j x_1^{(j)}t^j,\ldots, \sum_j x_N^{(j)}t^j)=\sum_j F^{(j)}t^j.$$

For the simplicity of the notation, we write $\bx^{(j)}=(x_1^{(j)},\ldots,x_N^{(j)})$.
For a point $P\in X$, let $P_m\in X_m$ be the trivial $m$-jet at $P$.
In particular if $P$ is the origin $0\in X\subset \bA_k^N$, then $0_m$ is defined by the maximal ideal $(\bx^{(0)},\ldots,\bx^{(m)})\subset k[\bx^{(0)},\ldots,\bx^{(m)}]$ in $\bA_k^{(m+1)N}$.
\end{say}

\begin{say}
\label{Frobenius}
 The Frobenius map of rings of positive characteristic has been important tool to study the singularities of positive characteristic. 
The concepts $F$-pure, strongly $F$-regular, weakly $F$-regular and $F$-rational appear in this stream. 
 These notions have close relations with rationality and log-canonicity:
A singularity is of
 dense $F$-rational type (i.e., it is $F$-rational by the reduction to characteristic $p$ for infinitely many prime number $p$) 
if and only if it is rational by Smith \cite{sm}, Hara \cite{ha}, Mehta and Srinivas \cite{ms}.
If a normal $\bQ$-Gorenstein singularity is of dense $F$-pure type
 (i.e., it is $F$-pure by the reduction to characteristic $p$ for infinitely many prime number $p$), then it is
 log-canonical by Hara and Watanabe \cite{hw}.
In the Gorenstein case, the three notions strongly $F$-regular, weakly $F$-regular and $F$-rational coincide.
When  we restrict ourselves in the case of   a complete intersection,
 we call it just $F$-regular.
 %
The definitions of $F$-pure and $F$-regular can be found in the papers above and we do not repeat them here.
\end{say}

\begin{lem}
\label{local}
  The $m$-jet scheme $X_m$ is $F$-pure (resp. strongly $F$-regular, rational) along the fiber $\pi_m^{-1}(P)$ if and only if $X_m$ is $F$-pure (resp. strongly $F$-regular, rational) at $P_m$.
\end{lem}

\begin{pf}
  Note that these conditions, $F$-pure,  strongly $F$-regular, rational, are
  open conditions. 
  Therefore, if $X_m$ has one of these conditions at $P_m$, then $X_m$ has that on  an open neighborhood $U\subset X_m$ of $P_m$.
 Remember that the multiplicative algebraic group  $\bA_k^1\setminus \{0\}$ acts on $X_m$ and the closure of the orbit of every point $y$ in $\pi_m^{-1}(P)$ contains $P_m$ (see, for example, \cite{cr}).
 This shows that on $X_m$   there is an isomorphism which sends $y$ into $U$. 
 Hence, $X_m$ has the condition at $y$.
\end{pf}


\begin{lem}[\textup{\cite[Lemma 3.9]{Ta}}] 
\label{takagi}
Let $(R,\fm)$ be a local ring at a closed point of a non-singular variety
              over an algebraically closed field of characteristic p and
$I \subset R$ an ideal. Fix any ideal $\fa \subset R$ and any real number $t
\ge 0$. Write $S=R/I$.
\begin{enumerate}
\item
The pair $(S, (\fa S)^t)$ is F-pure if and only if for all large $q=p^e \gg
0$, $\fa^{\lfloor t(q-1) \rfloor}(I^{[q]}:I) \not\subset \fm^{[q]}$.
\item
The pair $(S, (\fa S)^t)$ is strongly F-regular if and only if for every
element $g \in R \setminus I$, there exists $q=p^e>0$ such that $g\fa^{\lceil
tq \rceil}(I^{[q]}:I) \not\subset \fm^{[q]}$.
\end{enumerate}
\end{lem}


In the case of  a complete intersection, we can regard the following criteria of Fedder type as the definition of $F$-pure and $F$-regular.

\begin{cor}
\label{fedder}
If $(S, \fm)$ is a regular local ring of characteristic $p>0$,
$f_1,f_2,\ldots, f_r$ is an $S$-sequence and $f=\prod_{i=1}^r f_i$,
then the following are equivalent:
\begin{enumerate}
\item[(i)] $S/(f_1,\ldots, f_r) $ is $F$-pure (resp. $F$-regular)
\item[(ii)] $f^{p-1}\not\in \fm^{[p]}$ (resp. for any non-zero $g\in S$, there is $q=p^e> 0$ such that $gf^{q-1}\not\in \fm^{[q]}$).
\end{enumerate}
\end{cor}

\begin{pf}
The statement on $F$-purity is in \cite{fed}.   
 If $I=(f_1,\ldots ,f_r)$ is generated by $S$-regular sequence, then
$I^{[q]}: I = I^{[q]} + f^{q-1}S$ and our assertion on $F$-regularity follows from Lemma~
\ref{takagi}.

\end{pf}

To apply the criteria, we need to show that 
our jet schemes are complete intersections. 
The following is a characteristic free statement  and is a refinement of a special case of the statement obtained by Musta\cedilla{t}\v{a} \cite{must01} for characteristic zero.

\begin{lem}
\label{degdim}
  Let $X$ be a hypersurface of $\bA_k^N$ defined by a homogeneous polynomial $f$ of degree $d$. 
Assume $X$ has an isolated singularity at the origin $0\in X$.
Then, it follows:
\begin{enumerate}
\item[(i)] 
   If $d\geq N$, then $X_m$ is not irreducible for every $m\geq N-1$;
\item[(ii)]
  If $d\leq N-1$, then $X_m$ is irreducible, therefore a complete intersection, for every $m\in \bN$.   
\end{enumerate}
\end{lem}

\begin{pf}
  First of all, we note that for a hypersurface $X$ with the isolated singularity at $0$, the jet scheme $X_m$ is irreducible if and only if 
 \begin{equation}
 \label{dim}
  \dim \pi_m^{-1}(0)< (m+1)(N-1).
 \end{equation}
Indeed, as $\dim \pi_m^{-1}(X\setminus \{0\})= (m+1)(N-1)$,  ``only if" part is trivial. 
For the ``if" part, note that $X_m$ is defined by $m+1$ equations in 
$\bA_k^{(m+1)N}$.
Therefore, every irreducible component of $X_m$ has dimension greater than or equal to $(m+1)(N-1)$.
If we assume the inequality (\ref{dim}), then  $\pi_m^{-1}(0)$ 
does not provide with an irreducible component of $X_m$.

For the proof of (i), assume $d\geq N$. 
The  fiber $\pi_m^{-1}(0)$ is defined by 
$$F^{(0)}(\bo), F^{(1)}(\bo,\bx^{(1)}), \ldots, F^{(m)}(\bo,\bx^{(1)},\ldots, \bx^{(m)})$$
on $\bA_k^{mN}=\spec k[\bx^{(1)},\ldots, \bx^{(m)}]$ and the first $d$ polynomials are trivial because $F^{(j)}$ is homogeneous of degree $d$ and of weight $j$, 
therefore
every monomial in $F^{(j)}$ ($j<d$) has the factor $x_i^{(0)}$ for some $i$, where the weight of a monomial $\prod x_i^{(j)}$ is defined as $\sum j$.

Therefore, for $m\geq N-1$
   $$\dim \pi_m^{-1}(0)\geq mN-\max\{0, (m+1)-d\}\geq (m+1)(N-1).$$
 Hence, $X_m$ is not irreducible for  $m\geq N-1$.
 
 For  the proof of (ii), assume $d\leq N-1$.
 For $m$ such that $m\leq d-1$, as we see in the previous argument,
 $\pi_m^{-1}(0)=\bA_k^{mN}$. 
 As $m\leq N-2$, the inequality (\ref{dim}) holds, therefore $X_m$ is irreducible.
 For $m$ such that $m\geq d$, 
 we will show it by induction on $m$.
 Assume that $X_0=X,\ldots, X_{m-1}$ are irreducible. 
 We note that for $j\geq d$
 $$F^{(j)}(\bo,\bx^{(1)},\ldots,\bx^{(j)})=F^{(j-d)}(\bx^{(1)},\ldots, \bx^{(j-d+1)}),$$
 because $f$ is homogeneous of degree $d$.
 Since $\pi_m^{-1}(0)$ is defined by $F^{(j)}(\bo,\bx^{(1)},\ldots,\bx^{(j)})$ ($j=d,\ldots, m$) in $\spec k[\bx^{(1)},\ldots, \bx^{(m)}]$,
 we obtain
 $$\pi_m^{-1}(0)=\spec k[\bx^{(1)},\ldots, \bx^{(m)}]/\left(F^{(j-d)}(\bx^{(1)},\ldots, \bx^{(j-d+1)})\right)_{j-d=0,..,m-d}$$
 $$\simeq X_{m-d}\times\bA_k^{(d-1)N}.$$
 By this we have   $\dim\pi_m^{-1}(0)= (m-d+1)(N-1)+(d-1)N$ and it follows the inequality (\ref{dim}).
 Now we obtain the irreducibility of $X_m$ and in this case we have the codimension of $X_m$ equal to the number of the defining equation in $\bA_k^{(m+1)N}$.
\end{pf}

\section{Singularities of the jet schemes}

\begin{defn} Under the notation in \ref{notation}, let $k$ be a field of characteristic zero and $p$ a prime number.  
  Let $\fm$ be the maximal ideal $(\bx^{(0)},\bx^{(1)},\ldots, \bx^{(m)})\subset  k[\bx^{(0)},\bx^{(1)},\ldots, \bx^{(m)}]$.
  Take a  polynomial $F $  in the ring.
  A monomial $\bx\in k[\bx^{(0)},\bx^{(1)},\ldots, \bx^{(m)}]$ is called a 
  {\it  good monomial} for  $(F, p)$ if  $\bx\not\in\fm^{[p]}$ and $\bx\in F^{p-1}$  by modulo $p$ reduction.
   Here ``$\bx\in F^{p-1}$" means $\bx$ appears in $F^{p-1}$ with non-zero  
   coefficient.
\end{defn}

\begin{thm}
\label{main}
  Let $X$ be a  hypersurface in $\bA_k^N$ over a field $k$ of characteristic $0$ defined by a  polynomial of general type of 
  degree $d$,
  i.e., 
  $$f(x_1, \ldots, x_N)=\sum_{i_1\leq i_2\leq \ldots\leq i_d}\xi_{i_1,i_2,\ldots, i_d}x_{i_1}x_{i_2}\cdots x_{i_d},$$
  where $\{\xi_{i_1,i_2,\ldots, i_d}\}$ are algebraically independent over $\bQ$.
  If $d^2\leq N$, then the jet scheme $X_m$ is dense $F$-regular type for every 
  $m\in \bN$.
\end{thm}

\begin{pf}
  Fix $m\in \bN$. 
  Let $p$ be a prime number satisfying $p> m(d-1)+d$.
  By Lemma \ref{degdim}, we may assume that $X_m$ is a complete intersection.
  For the polynomial $f\in  k[x_1,\ldots, x_N]$, let 
   $F^{(j)}\in  k[\bx^{(0)},\bx^{(1)},\ldots, \bx^{(j)}]$ be  as in \ref{notation}
  and put $F=\prod_{j=0}^m F^{(j)}$.
  Let $g$  be any polynomial in $  k[\bx^{(0)},\bx^{(1)},\ldots, \bx^{(m)}]$.
  We will show that there exist $e>0$  and a monomial $M\in F^{p^e-1}$ with $gM\not\in \fm^{[p^e]}$ by modulo $p$.
  For any $e$, we can decompose 
  $$p^e-1=(p^e-p^{e-1})+(p^{e-1}-p^{e-2})+(p^{e-2}-1).$$
  Let $a=p^e-p^{e-1}, b=p^{e-1}-p^{e-2}, c=p^{e-2}-1$.
  Define three monomials $L_1, L_2, L_3 \in F=\prod_{j=0}^m F^{(j)}$.
  First pick up the term  $L_1(j)$ from $F^{(j)}$ as follows:
   $$
  \begin{array}{lll}
  L_1(0)& =& \xi_{1,2,\ldots, d}x_1^{(0)}\cdot x_2^{(0)}\cdots x_d^{(0)}\\
   L_1(1)& =& \xi_{d+1,d+2,\ldots, 2d}x_{d+1}^{(0)}\cdot x_{d+2}^{(0)}\cdots x_{2d-1}^ {(0)}x_{2d}^{(1)}\\
   L_1(2)& =&\xi_{2d+1,2d+2,\ldots, 3d} x_{2d+1}^{(0)}\cdot x_{2d+2}^{(0)}\cdots x_{3d-2}^ {(0)}x_{3d-1}^ {(1)}x_{3d}^{(1)}\\
  \cdots& &\cdots\\
  L_1(d-1)& =& \xi_{d^2-d+1,d^2-d+2,\ldots, d^2}x_{d^2-d+1}^{(0)}\cdot x_{d^2-d+2}^{(1)}\cdots x_{d^2-1}^ {(1)}x_{d^2}^{(1)}\\
  L_1(d)& =& \xi_{1,2,\ldots, d}x_{1}^{(1)}\cdot x_{2}^{(1)}\cdots x_{d-1}^ {(1)}x_{d}^{(1)}\\
  L_1(d+1)& =&\xi_{d+1,d+2,\ldots, 2d} x_{d+1}^{(1)}\cdot x_{d+2}^{(1)}\cdots x_{2d-1}^ {(1)}x_{2d}^{(2)}\\
  \cdots& & \cdots\\
  L_1(2d)& =& \xi_{1,2,\ldots, d}x_{1}^{(2)}\cdot x_{2}^{(2)}\cdots x_{2d-1}^ {(2)}x_{2d}^{(2)}\\
  \cdots& & \cdots\\
  \end{array}.$$
   Define $L_1=\prod_{j=0}^mL_1(j)/\mbox{(coefficients)}$.
  Then, note that every variable $x_i^{(j)}$ appears in $L_1$ at most once.
  We can see that $L_1^a\in F^a$ by modulo $p$, by noting that the coefficients 
  $\xi_{i_1,\ldots,i_d}$ are algebraically independent over $\q$.
 
   Next, pick up the term  $L_2(j)$ from $F^{(j)}$ as follows
$$ 
\begin{array}{lll} 
 L_2(0) &= & \xi_{2d,\ldots, 2d}(x_{2d}^{(0)})^d\\       
  L_2(1) &= &d \xi_{2d,\ldots, 2d}(x_{2d}^{(0)})^{d-1}x_{2d}^{(1)}\\ 
    L_2(2) &= &d  \xi_{2d,\ldots, 2d}(x_{2d}^{(0)})^{d-1}x_{2d}^{(2)}\\ 
   \cdots& & \cdots\\
    L_2(m) &=  & d\xi_{2d,\ldots, 2d}(x_{2d}^{(0)})^{d-1}x_{2d}^{(m)}\\  
 \end{array}.$$
    Define $L_2=\prod_{j=0}^mL_2(j)/\mbox{(coefficients)}$.
   Then, note that a variable with positive weight (i.e., $x_{2d}^{(j)}$ $(j>0)$) appears in $L_2$ at most once and 
   a variable $x_{2d}^{(0)}$ appears $m(d-1)+d$ times in  $L_2$.
   We can see that  $L_2^b\in F^b$ by modulo $p$.
  
  Finally, pick up the term  $L_3(j)$ from $F^{(j)}$  as follows:
  $$ 
\begin{array}{lll} 
 L_3(0) &=   & \xi_{d,\ldots, d}(x_d^{(0)})^d\\
  L_3(1) &=&d \xi_{d,\ldots, d}(x_d^{(0)})^{d-1}x_d^{(1)}\\   
    L_3(2) &=&d \xi_{d,\ldots, d}(x_d^{(0)})^{d-1}x_d^{(2)}\\          
   \cdots& & \cdots\\
   L_3(m) &=& d\xi_{d,\ldots, d}(x_d^{(0)})^{d-1}x_d^{(m)}\\           
    \end{array}.$$
   Define $L_3=\prod_{j=0}^mL_3(j)/\mbox{(coefficients)}$. 
  Then, note that a variable with positive weight (i.e., $x_{d}^{(j)}$ $(j>0)$) appears in $L_3$ at most once and 
   a variable $x_{d}^{(0)}$ appears $m(d-1)+d$ times in  $L_3$.
   We can see that  $L_3^c\in F^c$ by modulo $p$.
   
   Define $M:=L_1^a\cdot L_2^b\cdot L_3^c$.
   Noting that $(p^s-1)!$ has exactly $\left(\sum_{i=1}^{s-1}(p^i-1)\right)$-powers of $p$ as a factor and
   $(p^s-p^{s-1})!$ has exactly $(p^{s-1}-1)$-powers of $p$ as a factor for every positive integer $s$, we obtain that 
    $(p^e-1)!/(a!b!c!)$ does not have $p$ as a factor.
    Hence, it follows that $M\in F^{p^e-1}$ by modulo $p$.
   Every variable of weight $0$ appears in $M$ at most 
   $\max \{b(md-m+d), a+c(md-m+d)\}$ times and 
   $$p^e- \max \{ b(md-m+d), a+c(md-m+d)\} \to \infty \ \ (e\to \infty).$$
   On the other hand, every variable of positive weight appears in $M$ at most
   $a+b$ times (here, we used the fact that $a+c\leq a+b$).
   We can also see that 
   $$p^e-(a+b)\to \infty \ \ (e\to \infty).$$
   Therefore, for any polynomial $g\in  k[\bx^{(0)},\bx^{(1)},\ldots, \bx^{(m)}]$
   we obtain $gM\not\in \fm^{[p^e]}$ for sufficiently large $e$.
   \end{pf}

\begin{cor}[Theorem 1.2]
\label{rational}
 Let $k$ be a field of characteristic zero.
  Let $X$ be a  hypersurface in $\bA_k^N$ defined by a  homogeneous polynomial of general type of 
  degree $d$.
   If $d^2\leq N$, then the jet scheme $X_m$ has at worst rational singularities for every 
  $m\in \bN$.
\end{cor}

\begin{rem}
It is  expected that Theorem \ref{main}
and Corollary \ref{rational} also hold for the hypersurface $X$ of Fermat type of degree $d$
such that $d^2\leq N$.
\end{rem}

\begin{thm}
Assume $char k=p>0$.  Let $X$ be a hypersurface in $\bA_k^N$ defined by a  homogeneous polynomial $f\in k[x_1,\ldots, x_N]$ of degree $d$.
  If the jet scheme $X_m$ is a complete intersection and $F$-pure for every 
  $m\in \bN$, then $d^2\leq N$.
\end{thm}

\begin{pf}
As $X_m$ is a complete intersection and $F$-pure, 
there exists a good monomial $\bx$ for $(F,p)$.
Fix an expression  of $\bx$ into a product of monomials of $(F^{(j)})^{p-1}$'s.
Write $\bx=\prod_{j=0}^m\bx(j)$, where $\bx(j)$ is the contribution from $(F^{(j)})^{p-1}$.
Let $a_{ijk}$ be the power of $x^{(i)}_k$ in $\bx(j)$ and let $a_{ij}:=(\sum_{k=1}^N a_{ijk})/(p-1)$.
Then the matrix   $A=(a_{ij})_{0\leq i,j\leq m}$ satisfies the following conditions: 
\begin{enumerate}
\item
$A$ is an upper triangular matrix, 
\item
$ \sum_{i=0}^m a_{ij} =d, (  0\leq \forall j\leq m)$,
\item        
$ \sum_{i=0}^m  i a_{ij} = j , (0\leq \forall j\leq m)$ and
\item 
$ \sum_{j=0}^m   a_{ij} \leq N , (0\leq \forall i\leq m)$.
\end{enumerate}
  Under these conditions we will prove that  if $m$ is sufficiently large, for any real number $s < d^2$, there exists $i$
such that
   $ \alpha_i := \sum_{j=0}^m a_{ij} " s, $ which shows $d^2\leq N$ by (4).
  Let $C$ be the  matrix as follows:
     $$C=\left(
  \begin{array}{cccccccccccc}
  d & d-1&d-2 & \cdots & 1 &  & & & & & &  \\
     &  1  & 2 &               &d-1&d&d-1& \cdots &1& & & \\
     &      &    &               &     &  &  1 &            &d-1& d & d-1&\\
     &      &    &               &     &  &   &            &    &   & 1&\\
     &      &    &               &     &  &   &            &    &   &  & \ddots\\
 \end{array}
 \right)$$
In other words,    $C =  (c_{ij})_{0\le i,j \le m}$  be defined as 
 $ c_{ij} = d- u $   (if $ j=d i  \pm u$ )  for every $u=0,1,\ldots, d-1$ and
 $c_{ij}=0 $ (otherwise).
Then, $C$ also has the properties (1)$\sim$(3).
Let $ \gamma_i := \sum_{j=0}^m c_{ij}$,  and assume $m=dl$ for an integer $l>0$, then $\gamma_i=\frac{d(d+1)}{2}$
  (if $i=0$),
$d^2$  (if $1 \le i \le l-1$),
$\frac{d(d+1)}{2}$   (if $i =l$),
   0    (if $i \ge l+1$).

If we put
$\delta_i = \alpha_i - \gamma_i$  for  $i=0,\ldots ,m$,
then
  $$\sum_{i=0} ^m \delta_i  = 0,\ \ \operatorname{  and}\  \sum_{i=0} ^m  i \delta_i = 0.$$
 Now, assume $\max \alpha_i \le s = d^2 -e $ for some $e>0$.
Put $D = \sum_{i\le l - 1} \delta_i$ and $D' = \sum_{i\ge l } \delta_i$,
so that $D+D'=0$.
 Since $\delta_i \le \frac{d(d-1)}{2}-e$ for $i=0$ and $\delta_i \le -e$
for $1\le i \le l-1$, we have $D \le - le +\frac{d(d-1)}{2}$ and 
$D' \geq le -\frac{d(d-1)}{2}$.
  If we put $\delta_i= -e-\epsilon_i$
for $1\le i \le l-1$ and $D = - le +\frac{d(d-1)}{2}-\epsilon$, then $\epsilon\geq \sum_{i=1}^{l-1}\epsilon_i$ and $D' = le -\frac{d(d-1)}{2}+ \epsilon$.
By this we have
$$\sum_{i=0}^{l-1}i\delta_i=-e\sum_{i=0}^{l-1} i-\sum_{i=0}^{l-1}i\epsilon_i
\geq -e\frac{l(l-1)}{2}-(l-1)\epsilon.$$
  
On the other hand, noting that  $\delta_i\geq 0$ for $i\geq l+1$, we have 
$$\sum_{i=l}^mi\delta_i\geq lD'.$$  Thus, we conclude that
$$\sum_{i=0}^m i  \delta_i \ge - e\sum_{i=1}^{l-1}i + l(le -\frac{d(d-1)}{2})
 = \frac{1}{2}\left(el^2+(e-d^2+d))l\right).  $$
But if $l$ is sufficiently large, then the latter will be positive and
this contradicts to the fact that
$\sum_{i=0}^m i \delta_i=0$.
\end{pf}

\begin{say}
\label{fpt}
Takagi and Watanabe \cite{tw} introduced the invariant $F$-pure threshold (denoted by $fpt(X,Z,P)$) 
for a scheme $X$ over a field of positive characteristic and the closed subscheme $Z\subset X$ at a point $P\in X$.
It is closely related to the log-canonical threshold for characteristic zero.
Here, we refer the formula for a complete intersection case.

Let $k$ be a field of characteristic $p>0$.
Let $X$ be a subscheme of $\bA_k^N$ defined by polynomials 
$f_1,\ldots,f_r$ where $\dim X=N-r$.
Let $f=\prod_{i=1}^r f_i$. 
Let a closed subscheme $Z$ be defined by an ideal $I\subset k[x_1,\ldots, x_N]$ and $\fm$ be the maximal ideal of a point $P\in X$.
Let $q=p^e$.
Then, by Lemma \ref{takagi}
$$fpt(X,Z,P)=\lim_{q\to \infty}\frac{\max\{r\mid I^rf^{q-1}\not\subset \fm^{[q]}\}}{q}.
$$
As we think of only local a complete intersection case,
we can regard this formula as the definition of $F$-pure threshold.

\end{say}

\begin{thm}
  Let $X$ be a variety of locally a complete intersection at $0$ over a field of characteristic $p>0$.
  For $m<m'$, assume also $X_m, X_{m+1},\ldots, X_{m'}$ are  complete intersections at the trivial jets \\ $0_m,\ldots,  0_{m'}$. 
  Then, the following are equivalent:
\begin{enumerate}
\item[(i)]  
  $(X,0)$ is non-singular;
\item[(ii)]
   $fpt(X_m,0_m,0_m)=
  fpt(X_{m'},{\psi_{m'm}}^{-1}(0_m), 0_{m'})$, where $\psi_{m'm}:X_{m'}\to X_m$ is the truncation morphism.
\end{enumerate}
\end{thm}

\begin{pf}
   Assume (i), then $X_i$ is non-singular for every $i\in \bN$ and the 
   truncation morphism $\psi_{m'm}:X_{m'}\to X_m$ is smooth.
   In this case, $X_{m'}, X_m, {\psi_{m'm}}^{-1}(0_m)$ and $\{0_m\}$ are all non-singular. Therefore by the formula in \ref{fpt}, we have  $fpt(X_m,0_m,0_m)=\codim (\{0_m\}, X_m)=\codim ({\psi_{m'm}}^{-1}(0_m), X_{m'})= fpt(X_{m'},{\psi_{m'm}}^{-1}(0_m), 0_{m'})$.
  
   For the proof of (ii)$\Rightarrow$(i), 
   we first show 
\begin{equation}
\label{bigger}
  fpt(X_m,0_m,0_m)>
  fpt(X_{m+1},{\psi_{m+1,m}}^{-1}(0_m), 0_{m+1}),  
\end{equation}
if $(X,0)$ is singular and $X_m$ and $X_{m+1}$ are  complete intersections at the trivial jets.

   Let $\fm$ and $ \fm'$ be the maximal ideals of $\o_{X_m,0_m}$ and of $\o_{X_{m+1},0_{m+1}}$, respectively.
Let $f_1,\ldots, f_r$ define $X$ in $\bA^N$, where $r=\codim (X, \bA^N)$.
Then, under the notation in \ref{notation}, $X_i$ is defined by 
$F_l^{(j)}$ $(l=1,\ldots, r, j\leq i)$ in $\bA_k^{(i+1)N}=(\bA_k^N)_i$.
Let $G^{(j)}=\prod_{l=1}^r F_l^{(j)}$ and $G=\prod_{j=0}^m G^{(j)}$,
$G'=\prod_{j=0}^{m+1} G^{(j)}$.
For $q=p^e$,
 let 
 \begin{equation}
 \label{similar}
 r_q=\max\{s \mid \fm^sG^{q-1}\not\subset \fm^{[q]}\}
\end{equation}
\begin{equation}
\label{sim}
      r'_q=\max\{s \mid \fm^s{G'}^{q-1}\not\subset {\fm'}^{[q]}\}.
 \end{equation}     
Let $\bx$ be a monomial in ${G'}^{q-1}$ and $c$ be an element of 
$\fm^{r'_q}$ such that $c\bx\not\in {\fm'}^{[q]}$.      
Then $\bx$ is factored as $\bx=\bx'\bx''$, where $\bx'$ and $\bx''$ are   contributions from $G^{q-1}$ and  from $(G^{(m+1)})^{q-1}$, respectively. 
As $F_l^{(m+1)}$ is of weight $m+1$, each monomial of $F_l^{(m+1)}$ 
has at most one variable $x_i^{(m+1)}$ of weight $m+1$.
Then, if we factorize $\bx''=\bz\bz'$ with $\bz\in k[\bx^{(0)},\ldots,\bx^{(m)}]$ and $\bz'\in k[\bx^{(m+1)}]$, we have 
$$\deg \bz\geq (q-1)\sum_{j=1}^r(d_j-1),$$
where $d_j=\ord f_j$.
Here, we  note that $d_j\geq 1$ for all $j=1,\ldots, r$ and $d_j\geq 2$ for some $j$, since 
$(X,0)$ is singular.
The condition $c\bx\not\in {\fm'}^{[q]}$ gives $(c\bz)\bx'\not\in\fm^{[q]}$.
Noting that $\bx'\in G^{q-1} $ and $c\bz\in \fm^{r'_q+(q-1)\sum_{j=1}^r(d_j-1)}$,
we obtain 
$$r'_q+(q-1)\sum_{j=1}^r(d_j-1)\leq r_q,$$
which yields
$$fpt(X_{m+1},{\psi_{m+1,m}}^{-1}(0_m), 0_{m+1})+\sum_{j=1}^r(d_j-1)
 \leq fpt(X_m,0_m,0_m)$$
 as required in (\ref{bigger}).
 
 Now we can see the following in a similar and easier way as in the above discussions: 
 \begin{equation}
 \label{general}
 fpt(X_m, I,0_m)\geq fpt (X_{m+1}, I\o_{X_{m+1}}, 0_{m+1})
 \end{equation}
 for an ideal $I\subset \o_{X_m, 0_m}$.
 (This follows by just replacing $\fm^s$ by $I^s$ in (\ref{similar}) and (\ref{sim}).)
 By (\ref{bigger}) and (\ref{general}), we obtain that if $(X, 0)$ is singular, then $fpt(X_m,0_m,0_m)>
  fpt(X_{m'},{\psi_{m'm}}^{-1}(0_m), 0_{m'})$.
  This implies (ii)$\Rightarrow$(i).  \end{pf}

 \makeatletter \renewcommand{\@biblabel}[1]{\hfill#1.}\makeatother

\end{document}